\documentclass[12pt,reqno]{amsart} 

\usepackage{amssymb,bbm,enumitem,marginnote,url}
\usepackage{aliascnt}
\usepackage{doi}
\usepackage{pgfplots} 
\usepackage{booktabs}

\usepackage[margin=1.25in]{geometry}

\usepackage{psfrag}

\newcommand{\arxiv}[1]{%
 \href{https://arxiv.org/pdf/#1}{ArXiv:#1}}

\usepackage{color}

\newcommand\blue[1]{\textcolor{blue}{#1}}

\newtheorem{theorem}{Theorem}

\newaliascnt{lemma}{theorem}
\newtheorem{lemma}[lemma]{Lemma}
\aliascntresetthe{lemma}

\newaliascnt{proposition}{theorem}
\newtheorem{proposition}[proposition]{Proposition}
\aliascntresetthe{proposition}

\newaliascnt{corollary}{theorem}
\newtheorem{corollary}[corollary]{Corollary}
\aliascntresetthe{corollary}

\newaliascnt{conjecture}{theorem}

\aliascntresetthe{conjecture}

\newaliascnt{openQ}{theorem}

\aliascntresetthe{openQ}

\newaliascnt{quest}{theorem}

\aliascntresetthe{quest}

\newaliascnt{questx}{conjx}

\aliascntresetthe{questx}

\theoremstyle{definition}

\newaliascnt{defn}{theorem}
\newtheorem{defn}[defn]{Definition}
\aliascntresetthe{defn}

\newaliascnt{example}{theorem}

\aliascntresetthe{example}

\newaliascnt{rem}{theorem}

\aliascntresetthe{rem}

\makeatletter
\def\tagform@#1{\maketag@@@{\ignorespaces#1\unskip\@@italiccorr}}
\let\orgtheequation\theequation
\def\theequation{(\orgtheequation)}
\makeatother
\def\equationautorefname~{}

\newcommand{\B}{{\mathbb B}}

\newcommand{\R}{{\mathbb R}}

\newcommand{\N}{{\mathbb N}}
\newcommand{\C}{{\mathbb C}}
\newcommand{\Sp}{{\mathbb S}}

\newcommand{\V}{{\text{Vol}}}

\newcommand{\tr}{{\text{tr}}}
\newcommand{\Ric}{{\text{Ric}}}
\newcommand{\Div}{{\text{div}}}
\newcommand{\spec}{{\text{Spec}}}

\let\oldmarginnote\marginnote
\renewcommand{\marginnote}[1]{\oldmarginnote{\tiny \blue{#1}}}

\begin{document}
\title[Critical GJMS operators]{Extremal metrics of GJMS operators in critical dimensions}


\begin{abstract}
     Let $M$ be an even-dimensional closed Riemannian manifold. Graham-Jenne-Mason-Sparling constructed differential operators that are conformally covariant. When the dimension of the manifold coincides with twice of the order of the operator, it is called the critical GJMS operator. By generalizing the result of El Soufi and Ilias in dimension 2, and P\'{e}rez-Ayala in dimension 4, we are interested in finding necessary and partially sufficient conditions for extremal metrics for the eigenvalue functionals. Extremal metrics are  characterized by the stress energy tensor vanishing and the existence of the harmonic maps.
\end{abstract}
	
\author[]{Hanna N. Kim}
\address{Department of Mathematics, University of North Carolina at Chapel Hill, U.S.A.}
\email{hannakim@unc.edu}

	\maketitle 
\section{\bf GJMS operators in critical dimensions} \label{sec:intro}
In this section, we introduce the GJMS operators in critical dimensions. Let $M$ be a closed $n$-dimensional Riemannian manifold endowed with a metric $g$. Graham, Jenne, Mason and Sparling \cite{GJMS92} constructed formally self-adjoint, elliptic, conformally invariant differential operator of order $2N$ called GJMS operator denoted as $P_{N,g}$. These operators have conformal covariance structure and generalize this property for all dimensions. For example, $P_{1,g}$ is a Yamabe operator and $P_{2,g}$ is a Paneitz operator. In this paper, we mostly consider GJMS operators in critical dimensions ($2N=n$) due to their conformally invariant property with respect to the metrics. Moreover, it kills the constant functions. This is related to how $Q$-curvatures behave. Here, $Q$-curvature is a natural conformal covariant terms introduced by Branson \cite{B05}. 

By construction, $P_{N,g}$ acting on the smooth functions on $M$ has the following properties in the critical dimensions: For $u, w \in C^{\infty}(M)$,
\begin{enumerate}[label=(\roman*)]
\item $P_{N,g}$ is formally self-adjoint, that is, $\langle P_{N,g}\,u, w\rangle_{L^2(M,g)}=\langle \,u,P_{N,g}w\rangle_{L^2(M,g)}$.
\item $P_{N,g}$ is conformally invariant, that is, 
$$
P_{N,e^{2\omega}g}=e^{-2N\omega}P_{N,g}
$$
\item $P_{N,g}$ has polynomial expression in $\nabla$ and $R$ in which all coefficients are rational functions in $N$. Here, $\nabla$ is a covariant derivative and $R$ is a Riemann curvature tensor.
\item $P_{N,g}=\Delta^{N}+\text{l.o.t}$ where l.o.t means ``lower order terms". 
\item $P_{N,g}$ has the form $\delta S_{N,g}d$ where $S_{N,g}$ is an operator on 1-forms given as $\Delta^{N-1}+\text{l.o.t}$ .
\end{enumerate}
We can also express the operator as $(d\delta)^{N}+\text{l.o.t}$. For simplicity, we may write $P_{N,g}\equiv P_g$ from now on. By the conformal invariance property, for all $g_\omega \in [g]$ 
$$
\langle P_{g_\omega}u,w\rangle_{L^2(M,g_\omega)}=\langle P_{g}u,w\rangle_{L^2(M,g)}.
$$

For a closed manifold $M$, we have a sequence of increasing eigenvalues for $P_{g}$
    \[
        \lambda_1(P_g) \leq \lambda_2(P_g) \leq \cdots \leq\lambda_k(P_g)\leq \cdots \to \infty
    \]   
where each eigenvalue is repeated according to multiplicity. It is possible for $P_g$ to have negative eigenvalues. By the work of Canzani, Gover, Jakobson and Ponge \cite[Theorem 4.3]{CGJP14}, the number of negative eigenvalues is conformally invariant. As a remark, the kernel of $P_{g}$ containing nowhere vanishing function is equivalent to $P_g$ having constant function in its kernel. 

We want to study the properties of critical points of the $k$-th eigenvalue functional associated with $P_g$ written as 
    \begin{equation} \label{evfunctional}
        g\mapsto \lambda_k(P_{g}).
    \end{equation}
Here, we assume that the volume is preserved under the perturbation. In particular, we can ask this question inside the conformal class. For any background metric $g$, we write $g_\omega \in [g]$ where $\omega$ is a conformal factor. And we are interested in the functional
    \begin{equation} \label{evfunctionalconf}
    e^{2\omega}g\in [g]\mapsto \lambda_k(P_{e^{2\omega}g}).
    \end{equation}
    where we set $\V(M,g_\omega)=\V(M,g)$. Before we state our main goal, let's briefly introduce how similar terms have been defined. 
\begin{itemize}
    \item \textit{C-extremal}: El Soufi and Ilias \cite{EI03} defined $C$-$extremal$ metric in a class of metrics $C(g)$ that is conformal to a metric $g$ with normalized volume if the first eigenvalue of Laplace-Beltrami operator is a local maximizer for any perturbation of metrics.
        
    \item \textit{conformally extremal for the $k$-th normalized eigenvalue functional or conformally extremal for $\lambda_k$}: P\'{e}rez-Ayala \cite{P22} defined normalized $k$-th eigenvalue functional associated with the Paneitz operator $P_{4,g}$ in dimension 4.

    \item \textit{extremals of renormalised eigenvalue}: Humbert, Petrides and Premoselli \cite{HPP25} studied the case of eigenvalue optimization when the order of the GJMS operator of order $N$ when the dimension of the manifold is $n$ with $2N<n$.
\end{itemize}

As a remark, the term \textit{conformally maximal} was used in the work by \cite{KNPP19} as a global maximizer in a conformal class. We use the following definition for our case:

\begin{defn} \label{confext}
\begin{enumerate}[label=(\roman*), leftmargin=1.5em, labelsep=0.4em, nosep]
    \item If the $k$-th eigenvalue functional in \autoref{evfunctional} has a critical point at some metric, we say that the metric is \textit{extremal} for the $k$-th eigenvalue (normalized) functional. 
    \item If the $k$-th eigenvalue functional in \autoref{evfunctionalconf} has a critical point at some metric in the conformal class, we say that the metric is \textit{conformally extremal} for the $k$-th eigenvalue (normalized) functional. 
\end{enumerate}
\end{defn}

We define what it means to be critical point at some metric later. In this paper, we discuss the extremal metrics of GJMS operators in critical dimensions. El Soufi and Ilias \cite{EI08} in 2-dimensional and P\'{e}rez-Ayala \cite{P22} in 4-dimensional closed manifolds, studied the extremal metrics of the Laplace-Beltrami and the Paneitz operator, respectively. 
We generalize their results to higher order of GJMS operators in critical dimensions by relying on the fact that the coefficients of $P_g$ can be written as a universal polynomial of $g_{ij}$, $g^{ij}$ and its partial derivatives. 

For the functional \autoref{evfunctional}, there exists a traceless 2-tensor associated to the corresponding eigenfunctions. When the metric is critical, the sum of these tensors vanish and there exists a harmonic map consisting of the eigenfunctions as its coordinates. 

\section{Main results} \label{wodz}

This section first discusses characterizing conformally extremal metrics in critical GJMS operators and then characterizing these metrics. The difference is that for the first case, the perturbation stays within the conformal class and for the second case, the perturbation varies within all smooth Riemannian metrics. We use the fact that GJMS operators of orders $N\leq n/2$ of $n$-dimensional closed manifolds can be written as linear operators with coefficients depending analytically on the metric. 

\subsection{Conformally extremal metrics for critical GJMS operators}

We introduce a bilinear differential operator acting on $C^{\infty}(M)\times C^{\infty}(M)$ of order $2N$ that describes GJMS operators in critical dimensions. 
\begin{defn} \label{bilinear}
    Let $(M,g)$ be a $2N$-dimensional and closed manifold with a Riemannian metric $g$ and $P_{N,g}$ be a critical GJMS operator on $M$. Define a bilinear differential operator $B_{N,g}:C^{\infty}(M)\times C^{\infty}(M) \to C^{\infty}(M) $ as
       \begin{equation}
       2B_{N,g}(u,w)=uP_{N,g}(w)+wP_{N,g}(u)-P_{N,g}(uw),\label{PandB} 
    \end{equation}
    where $u,w \in C^\infty(M)$.
\end{defn}
This is essentially a Leibniz product rule of $P_{N,g}$ where the remainder is $B_{N,g}$. We write $B_{N,g}=B_g$ for simplicity. This bilinear differential operator is symmetric and conformally covariant ($B_{e^{2\omega}g}=e^{-2N\omega}B_g$). The symmetry is due to the symmetry of $P_{g}$ and the fact that $1 \in \ker(P_g)$ for critical GJMS operators. Indeed, we have
    $$
    \int_M B_g(u,w)dv_g= \int_M uP_g(w)dv_g.
    $$
As a remark, it is not true in general that $1$ belongs in the kernel when GJMS operator is not critical.

When the metric is conformally flat, we have more information for $B_g$. In dimensions 4, Connes \cite{C95} used quantized calculus to compute a natural bilinear operator for the Paneitz operator. Ugalde \cite{U06} generalized his result and constructs a bilinear operator based on Wodzicki's residue of a commutator operator acting on middle dimension forms for all critical GJMS operators. 

In general, the explicit expression of $B_{g}$ becomes complicated. Juhl \cite{J13} proposed recursive formulas for computing GJMS operators with lower order GJMS operators. An explicit formula for $N=3$ was constructed by Branson \cite{B85} and W\"{u}nsch \cite{W86}. Also, see the work by Chen and Hou \cite{CH16}. 

The next results in this paper are the necessary and partial sufficient conditions of conformally extremal metrics. The conditions rely on the fact that the relationship between conformally extremal metrics and the set of eigenfunctions that maps the manifold to a higher dimensional sphere. 

Here, the conditions are generalization of the Laplace Beltrami operator in the paper by El Soufi and Ilias \cite{EI03} and Paneitz operator by P\'{e}rez-Ayala \cite{P22}. We note that Petrides and Tewodrose \cite[Theorem 6.2]{PT25} developed a framework of similar understanding of critical metrics using Clarke’s differentiation theory for all conformal GJMS operators. They also showed the following result in their paper.

Let's discuss the necessary conditions of conformally extremal metrics first.
\begin{theorem}[Necessary condition]    \label{main}
    Let $(M, g)$ be a $2N$-dimensional closed Riemannian manifold with a background metric $g$. Suppose that $g_e \in [g]$ is conformally extremal for nonzero $\lambda_k(P_{g_e})$. Then there exists a finite family of eigenfunctions $\{u_1, \dots u_\ell\}$ associated with $\lambda_k(g_e)$ satisfying 
    \begin{equation}
        \sum_{i=1}^\ell u_i^2\equiv 1. \label{constantcon}
    \end{equation}
\end{theorem}
The property \autoref{constantcon} plays an important role in the eigenvalue optimization. Li and Yau \cite{LY82} in 2 dimensions and El Soufi and Ilias \cite{EI86} in all dimensions but not less than 2, showed that in symmetric spaces, the first non-zero Laplace-Beltrami eigenvalue can be maximized with the induced metric with a harmonic minimal immersion. 

Two immediate necessary conditions following this theorem are that the eigenvalue corresponding to conformally extremal metric has at least multiplicity 2, and we can characterize the eigenvalue with the sum of the associated bilinear map $B_g$. 

Here, we say $\lambda_k(P_g)$ is \textit{degenerate} if the eigenspace of $\lambda_k(P_g)$ denoted as $E_k(g)$ has dimension greater than 1. This means that either $\lambda_k(P_g)=\lambda_{k+1}(P_g)$ or $\lambda_{k-1}(P_g)=\lambda_{k}(P_g)$ holds.

\begin{corollary} \label{cormain}
    Let $(M, g)$ be a $2N$-dimensional closed Riemannian manifold with a background metric $g$. Suppose that $g_e \in [g]$ is conformally extremal for $\lambda_k(P_g)$ where $\lambda_k(P_{g_e})\neq 0$. Then $\lambda_k(P_{g_e})$ is degenerate and 
    \begin{equation}
        \lambda_k(P_{g_e})=\sum_{i=1}^\ell B_g(u_i,u_i) \label{eigformula}
    \end{equation}
    where $B_g$ is the bilinear differential operator defined in \autoref{bilinear}.
\end{corollary}

Next, we introduce a density function which is a direct application of the previous result. See the work by Angelsberg \cite{A07}, B\'{e}rald \cite{B08}, Biquard and Madani \cite{BM12}, Lin and Zhu \cite{LZ25} and P\'{e}rez-Ayala \cite[Corollary 1.2]{P22} for the Paneitz operator. In our case, the eigenfunctions naturally form a conformal-harmonic map to a sphere of dimension $\ell-1$. This gives another necessary condition, which is equivalent to \autoref{main}.

Consider a sphere-valued map $U=(u_1,\dots,u_\ell):(M,g) \to \Sp^{\ell-1}$ and suppose that this map $U$ is a solution to the equation
$$
P_g(U)=\mathcal{E}_g(U)U
$$
where $\mathcal{E}_g$ is an energy density function defined by the bilinear differential operator as 
\begin{equation}
    \mathcal{E}_g(U)=\sum_{i=1}^{\ell}B_g(u_i,u_i). \label{energy}
\end{equation}
We call such a map $U$ \textit{conformal-harmonic map}. Indeed, suppose that there exists a metric $g_e$ satisfying the equation \autoref{main} under the same assumption. Then $U$ consisting of eigenfunctions maps $M$ to a sphere $\Sp^{\ell-1}$ and the energy density $\mathcal{E}_g$ defined by the bilinear operator is conformally covariant since $B_g$ is conformally covariant.  

\begin{corollary}
    Under the same assumptions of \autoref{main}, the map $U=(u_1, \cdots, u_\ell): (M, g_e) \to \Sp^{\ell-1}$ is a conformal-harmonic map.
\end{corollary}

Next, the partial sufficient condition says that if there exists a metric with nonzero eigenvalue with additional condition, the necessary condition becomes a sufficient condition as well.
 
For the results regarding partial sufficient conditions in dimensions 2, see El Soufi and Ilias \cite{EI03}[Theorem 2.1], in dimensions 4, see P\'{e}rez-Ayala \cite[Theorem 1.3 (ii)]{P22}. Karpukhin and M\'{e}tras \cite{KM22} generalized this map of Laplace Beltrami operator in all dimensions. 
    \begin{theorem}[Partial sufficient condition] \label{partialsuf}
    Let $(M,g)$ be a $2N$-dimensional closed Riemannian manifold with a background metric $g$. Assume there is a metric $g_e \in [g]$ such that $\lambda_k(P_{g_e}) \neq 0$ and for which either $\lambda_k(P_{g_e})>\lambda_{k-1}(P_{g_e})$ or $\lambda_k(P_{g_e})<\lambda_{k+1}(P_{g_e})$ holds. If there exists a family of eigenfunctions $\{u_1, \dots, u_\ell\}$ associated with $\lambda_k(P_{g_e})$ such that $\Sigma_{i=1}^\ell u_i^2$ is constant, then $g_e$ is conformally extremal for $\lambda_k$.
    \end{theorem}
    Here, note that the gap condition of the $k$-th eigenvalue with its neighboring eigenvalue in \autoref{partialsuf} is necessary since we need to consider left and right derivatives of the eigenvalues. See \autoref{in} for details later. Using this theorem, we can also understand this using the energy density that we defined in \autoref{energy}.
    \begin{corollary}\label{partialsufcor}
        Let $(M,g)$ be a $2N$-dimensional closed Riemannian manifold with a background metric $g$. Assume that there is a smooth conformal-harmonic map $U:(M,g) \to (\mathbb{S}^{\ell-1},g_r)$ with energy density satisfying $\mathcal{E}_g(U)\neq 0$ pointwise on $M$. Then there exists a smooth metric $g_U$ in $[g]$ such that $g_U$ is conformally extremal for $\lambda_k$ for some $k$, and $\mathcal{E}_{g_U}(U)=\lambda_k(P_{g_U})$.
    \end{corollary}

In the work by Ammann and Jammes \cite{AJ11}, they showed that GJMS operator (in fact, for conformally covariant elliptic operator) whose order is strictly smaller than the half of the dimension of the manifold, the supremum among the conformal class is always $+\infty$. However, this case does not apply to our case since the order $N$ is equal to the half of the dimension $2N$. We have the following upper bound analogous to the results of Korevaar \cite{K93} and Hassannezhad \cite{H11}, and direct generalization of P\'{e}rez-Ayala \cite{P22}.
    \begin{theorem} \label{main2}
        Let $M$ be a $2N$-dimensional closed Riemannian manifold with a conformal class $[g]$. Then there exists a constant $A=A([g])>0$, depending only on the conformal class of $g$, such that
        $$
            \sup_{g_w \in [g]}\lambda_k(P_{g_w})\mathrm{Vol}(M,g_\omega) \leq A(k+1)< \infty     
        $$
    \end{theorem}
This result provides a background on that you can find the global maximum in the conformal class. The homogeneous manifolds were studied in El Soufi and Ilias \cite{EI08}[Corollary 2.1] for Laplace Beltrami operator and the case of sphere was studied for Paneitz operator in the work by P\'{e}rez-Ayala \cite[Theorem 1.6]{P22}.

Using the result above, we maximize the first nonzero eigenvalue for the even-dimensional spheres. It's important to know the sign of the eigenvalue. Kazdan and Warner \cite{KW75} showed that the sign of the first eigenvalue $\lambda_1$ is conformally invariant for Yamabe operator. Canzani, Gover, Jakobson and Ponge \cite{CGJP14} generalized this result to GJMS operators. It is known that there exists a large class of conformal classes with first positive eigenvalue in dimensions 4 by the work of Gursky and Viaclovsky \cite{GV03}. However, for critical GJMS operator of order 6 and higher, it is still not well-known in general manifolds. 

In our case, we use the following fact that if $g$ is an Einstein metric with $\Ric_g=\alpha(2N-1)g$, one can write 
\begin{equation}\label{einstein}
    \displaystyle {P_{N,g}=\Pi_{1\leq j\leq N}\left(\Delta_g+\alpha(N+j-1)(N-j)\right)}
\end{equation}
where this calculation was shown by Graham \cite{G03}, Gover \cite{G06} (See the paper by Guillarmou-Naud \cite{GN06} for constant sectional curvature spaces). For a conformal class of round metric, note that $\text{Ker}(P_{N,g})$ consists of only constant functions. So it makes sense to find the upper bound of the second eigenvalue. 

The first nonzero eigenvalue can be shown to have its maximizer using the Hersch's idea in 2 dimensions \cite{H70} and El Soufi and Ilias's idea in all higher dimensions \cite{EI08}. 
\begin{theorem}
    \label{first}
    For $2N$-dimensional sphere $\Sp^{2N}$, consider a conformal class of round metric denoted as $[g_0]$. If we normalize the metric $g_\omega \in [g_0]$ so that $\V(\Sp^{2N},g_\omega)=\V(\Sp^{2N},g_0)$, then we have 
    \begin{equation*}
        \lambda_2(P_{g_\omega}) \leq (2N)!\label{SecEig}
    \end{equation*}
    with equality if $g_\omega$ is isometric to $g_0$. 
\end{theorem}
The upper bound for the next higher eigenvalue has been studied the folding method (See the work by Girouard, Nadirashvili and Polterovich \cite{GNP21}, Kim \cite{K25}) in case of the Laplace-Beltrami operator. As a remark, one can ask questions similar questions for symmetric manifolds in even dimensions for Laplace-Beltrami operator, which was studied by Eddaoudi and Girouard \cite{EG25}.

    \subsection{Critical metrics for critical GJMS operators}
    Our main result characterizes the critical metrics among analytic perturbation. 

    For the Laplace-Beltrami operator, El Soufi and Ilias \cite{EI08} obtained a the necessary condition 
    \begin{equation*} \label{takahashi}
        \sum_{i=1}^\ell du \otimes du=g,
    \end{equation*}
which says that when finding critical metrics in all Riemannian metrics, one needs to find eigenfunctions that admits isometric minimal immersion to a sphere (See Takahashi \cite{T66}). The tensor defined as $du \otimes du-\frac{1}{2} |\nabla_g u|^2g$ is known as the stress-energy tensor. We derive a similar tensor condition for the necessary condition and sufficient condition. 
	\begin{theorem}\label{necessary2}
   	 Let $g$ be a Riemannian metric for which $\lambda_k(P_g)\neq 0$ and for which either $\lambda_k(P_g)>\lambda_{k-1}(P_g)$ or $\lambda_k(P_g)<\lambda_{k+1}(P_g)$ holds. 
     \begin{enumerate}[label=(\roman*)]
         \item If the metric $g$ is critical for the normalized $k$-th eigenvalue functional, then there exists a family of eigenfunctions $\{u_1,\cdots,u_\ell\}$ associated with $\lambda_k(P_g)$ such that \begin{equation}\label{summation}
         \sum_{i=1}^\ell\tau_g(u_i,u_i)=0 \qquad
		\sum_{i=1}^\ell u_i^2=1,
		\end{equation}
            where $\tau_g$ is a symmetric traceless 2-tensor determined by $P_g$.
         \item If there exists a set of eigenfunctions $\{u_1, \dots, u_\ell\}$ associated with $\lambda_k(P_g)$
         that satisfies \autoref{summation}, then $g$ is critical for the $k$-th eigenvalue functional.     \end{enumerate}		
	\end{theorem}


\section{Perturbation theory and critical metrics}

Here, we consider GJMS operators of orders $N=1,\dots,n/2 $ for $n$ dimensional manifolds when $n$ is even. In this section, we study the critical metrics of the GJMS operator. The existence of critical metrics was studied by El Soufi and Ilias \cite{EI08} for the Laplace Beltrami operator and by P\'{e}rez-Ayala \cite{P22} for the Paneitz operator. In general, the functions $t\mapsto\lambda_k(P_{g(t)})$ is continuous but not differentiable. Canzani \cite{C14} established that the eigenvalue of GJMS operator is analytic under conformal perturbation of metrics based on the work of Rellich \cite[page 74]{R69}. Here, we consider an analytic perturbation of metrics that is not necessarily in the conformal class.

By the perturbation theory, we define left and right derivatives of $t \mapsto \lambda_k(P_{g(t)})$ along some locally analytic curve of metric deformation. We write analytic metric deformation of the metric as $g(t)$ such that $g(0)=g$. The following theorem establishes well-definedness of the mapping from $t$ to metric, to its corresponding operator, eigenvalues and eigenfunctions for small time.

\begin{proposition}\label{onesided}
    Denote $P_{g(t)}$ as a family of GJMS operators of order $N$ where $g(t)$ is an analytic family of the metric with $g(0)=g_0$. Suppose that $\lambda_k(P_{g_0})\equiv\lambda_k$ is an eigenvalue with multiplicity $\ell$ of the operator $P_{g_0}$ and  there is a positive number $\tilde{\epsilon}$ such that $(\lambda_k-\tilde{\epsilon},\lambda_k+\tilde{\epsilon})$ consists of exactly $\lambda_k$. Then there exist analytic real-valued functions $\Lambda_1(t), \dots ,\Lambda_\ell(t)$ and $u_1(t),\dots u_\ell(t)$ for small $|t|$ satisfying 
    \begin{enumerate} [label=(\roman*)]
        \item For $i=1,\dots \ell$, $u_i(t)$ is an eigenfunction of $P_{g(t)}$ with an eigenvalue $\Lambda_i(t)$ where $\Lambda_i(0)=\lambda_k$. 
        \item $\{u_1(t),\dots,u_\ell(t)\}$ is an orthonormal set for each $t$.
        \item For each interval $[\lambda_k-\epsilon, \lambda_k+\epsilon]$ where $\epsilon<\tilde \epsilon$, there exists a number $\delta >0$ such that the intersection of the spectrum of $P_{g(t)}$ with this interval has only $\Lambda_1(t),\dots,\Lambda_{\ell}(t)$, for $|t|< \delta$.
    \end{enumerate}
\end{proposition}
 By $(i)$, we can write 
\begin{equation}\label{varlaplace} 
    P_{g(t)}u_i(t)=\Lambda_i(t)u_i(t), \qquad \Lambda_i(0)=\lambda_k(P_{g})
\end{equation}
for $i=1,\cdots,\ell$ and the orthonormal property in $(ii)$
    \begin{equation}\label{normalize}
    \int_M u_i(t)u_j(t)dv_{g(t)}=\delta_{ij}
    \end{equation}
for $i,j \in\{1,\cdots,\ell\}$.
Since this theorem allows us to locally perturb the eigenvalues $\Lambda_i(t)$, we can separate the spectrum along the perturbation, which is $(iii)$.

    The proof of $(i)$ and $(ii)$ in \autoref{onesided} proceeds similarly to the work of Bando and Urakawa \cite{BU83}[Theorem A.3]. To show this, we first verify that $P_{g(t)}$ is analytic based on the recursive formula developed by \cite{J13}. Then we use analyticity of the coefficients of $P_{g(t)}$ by \cite{BU83}[Lemma A.4] to show that $P_{g(t)}$ is analytic. Then using \cite{BU83}[Theorem A.2], we have $(i)$ and $(ii)$ in \autoref{onesided} by verifying the conditions. For $(iii)$, we adapt the method by Canzani \cite[Proposition 7.3]{C14} that there exists $\epsilon$ such that $|t| <\delta$ gives that 
    $$
    \text{Spec}(P_{g(t)}) \cap [\lambda_k-\epsilon,\lambda_k+\epsilon]=\{\Lambda_1(t),\dots,\Lambda_\ell(t)\}
    $$
    for all $|t|<\delta$ and where the eigenvalues $\Lambda_i(t)$'s are analytic functions. Here, the functions $\Lambda_i(t)$'s can be different depending on the choice of each branch. 

    We introduce some criterion first. Let $X$, $Y$ be complex Banach spaces. We say that a family of bounded operators $A(t)$, $|t|< \delta$ for small $ \delta$, of $X$ into $Y$ is analytic if it can be expanded as convergent power series of bounded operators $A^{(1)}$, $A^{(2)}$, $\dots$ of $X$ into $Y$ such that
    \begin{equation*}
        A(t)=\sum_{k\ge0}A^{(k)}t^k
    \end{equation*}
    for every $|t|< \delta $. The series converges in uniform topology.
    
    For a function $u$ on $M$ for $t$ with $|t|<\delta$, write the multiplication operator $m(t):H^{2N}(M,g(t))\to L^2(M,g_0)$ 
        \begin{equation}
            m(t)u=(\det(g(t))/\det(g_0))^{1/4}u \label{multiop}.
        \end{equation}
    By definition, they are bounded isomorphism of the Hilbert space $H^{2N}(M)$ to itself. Define the operator $A(t): H^{2N}(M,g(t))\to H^{0}(M,g(t))$ as
        \begin{equation}
        A(t)=m(t)P_{g(t)}m(t)^{-1} \label{Aop}
        \end{equation}
        where $A(0)=P_{g}$. Indeed, $A(t)$ is a self-adjoint and symmetric operator. Moreover, $m(t)$ sends the eigenfunctions of $g(t)$ to the eigenfunctions of $A(t)$ with the same eigenvalue.
        
    We prove the following proposition later in \autoref{computation}.
    \begin{proposition} \label{Panalytic}
        Let $(M,g_0)$ be an $n$-dimensional manifold with a Riemannian metric $g_0$ for even $n$. Given an analytic perturbation of metrics $g(t)$ where $g(0)=g$, Write GJMS operators $P_{2N,g(t)}$ depending on this perturbation where $N=1,\dots, n/2$. Then $P_{2N,g(t)}$ is analytic on $|t|<\delta$ for some $\delta$.
    \end{proposition}

    Based on this, we can define the analytic perturbation of the metrics and operators. For bounded operators, the operator converges in the uniform topology. 

    We adapt the proof by Canzani \cite[Proposition 7.3]{C14}. For a fixed $z \in \C$, let $A_t(z):=A(t)-z$ where $A(t)$ is defined in \autoref{Aop}. It is well known that if $z$ is in the resolvent of $A(t)$ ($z$ is not an eigenvalue of $A(t)$), then we know that $A_t(z)$ is surjective. Moreover, for $z_0$ in the resolvent set of $P_{g_0}$, we have 
    \begin{equation} \label{b0bound}
        ||(P_{g_0}-z_0)u|| \ge b_{g_0}||u||\qquad \text{where $b_{g_0}=\inf_{\lambda\in\spec (P_{g_0})}|\lambda-z_0|$}.
    \end{equation}
    where $||\cdot||$ is a $L^2$-norm.
    \begin{lemma} \label{resolvent}
        There exists $\gamma,\delta>0$ where the resolvent operator $R_t(z):=A_t(z)^{-1}$ exists for $|t|<\delta$ and $|z-z_0|<\gamma$. Moreover, for every $u \in C^\infty(M,g_0)$, $R_t(z)u$ depends continuously on $z$ and $t$ in the $L^2$-norm.
    \end{lemma}
         Note that $A(t)$ is also an elliptic operator since $P_{g(t)}$ is elliptic. For all $u \in C^\infty(M)$, there exists a constant $c_k$ independent of $u$ such that 
    \begin{equation} \label{sobolbound}
        ||u||_{k+m}^2\leq c_k(||P_{g_0}u||_k^2+||u||_0^2)
    \end{equation}
   (See \cite[p.193]{LM89}). Here, $||\cdot||_k$ is equivalent to the $k$-Sobolev norm.
    \begin{proof}
    It suffices to show the injectivity of $A(t)$. For all $\eta>0$, there exists $\gamma,\delta>0$ such that
    $$
    ||A_t(z)u||\ge (b_{g_0}-\eta)||u||,
    $$
    where $|z-z_0|<\gamma$ and $|t|<\delta$. We prove by contradiction. Suppose that there exists a sequence of sets $(\gamma_i,\delta_i,u_i, z_i)$ where $\gamma_i, \delta_i \to 0$, 
    $$
    ||(A_{t_i}(z_i)u_i||< (b_{g_0}-\eta)||u_i||    
    $$
    for some $\eta>0$.
    Using \autoref{sobolbound} and \autoref{b0bound}, we write
    \begin{align*}
        ||((P_{g_0}-z_0)-A_{t_i}(z_i))u_i||
        &\ge  ||((P_{g_0}-z_0)u_i||-||A_{t_i}(z_i)u_i||\\ 
        &\ge b_{g_0}-(b_{g_0}-\eta)=\eta,
    \end{align*} 
    where the last line is due to normalizing $u_i$ so that $||u_i||=1$. This contradicts that $||((P_{g_0}-z_0)-A_{t_i}(z_i))u_i|| \to 0$ since $A(t)$ is analytic and $u$ is bounded by \autoref{sobolbound}. Then the surjectivity for $|t|<\delta$ follows from this. Indeed, suppose not, then  $z \in \spec(P_{g(t)})$. From $A_t(z)u=0$ for the corresponding eigenfunction $u$, it follows that $u=0$, which is a contradiction. 

    Define the resolvent as $R_t(z)=(A_t-z)^{-1}$. To show the continuity of the $R_t(z)u$ with respect to $z$ and $t$, we write
    \begin{align}
        ||R_t(z)u-R_0(z_0)u||&\leq \frac{1}{b_{g_0}-\eta}||R_t^{-1}(z)R_t(z)u-R_t^{-1}(z)R_{0}(z_0)u||\nonumber\\
        &=\frac{1}{b_{g_0}-\eta}||R_0^{-1}(z_0)R_0(z_0)u-R_t^{-1}(z_0)R_{0}(z_0)u||\nonumber\\
        &\le\frac{1}{b_{g_0}-\eta}(||P_{g_0}-A(t)||||R_0(z_0)u||+||z-z_0||||R_0(z_0)u||) \nonumber
    \end{align}
    and the continuity on $z$ and $t$ follows since $A(t)$ is continuous on $g$.
    \end{proof}
     For a fixed $g_0$, let $g(t)$ be an analytic perturbation of $g_0$ for all $|t|<\delta$. Recall the multiplication operator defined in \autoref{multiop}. Choose $\delta$ so that the isomorphism $m(t):C^\infty(M,g(t)) \to C^\infty(M,g_0)$ is defined for all $|t|<\delta$. Let $\Gamma$ be a differentiable curve with interior domain $D \subset \C$. For fixed $t$, write $\mathcal{F}_t$  which is a linear subspace of $C^\infty(M,g_0)$ as
    $$
    \mathcal{F}_t=\text{span}\{ m(t) u:u\in \ker(P_g-\lambda I) \text{ for $\lambda \in D \cap \spec (P_{g(t)})$} \}
    $$
    \begin{lemma} \label{separate}
        If $\Gamma$ intersects with none of the eigenvalues of $P_{g_0}$, then there exists $\delta$ so that for all $g(t)$ where $|t|<\delta$, there is exactly $\ell$ number of
 eigenvalues of $P_{g(t)}$ in the interval $(\lambda-\epsilon, \lambda+\epsilon)$.
    \end{lemma}
    \begin{proof}
        For $|t|< \delta$, define the spectral projection operator $F_t(\Gamma): C^\infty(M,g_0) \to C^\infty(M,g_0)$ to be the projection of $C^\infty(M,g_0)$ onto $\mathcal{F}_{t}(\Gamma)$. Since  we assume $\Gamma \cap \spec (P_{g_0}) =\emptyset$, by \autoref{resolvent}, there exists a neighboring curve $\Gamma'$ near $\Gamma$ and $\delta<\tilde \delta$ so that $\spec (P_{g(t)}) \cap \Gamma'$ is empty.  By holomorphic functional calculus, we have 
        $$
        F_{t}(\Gamma)u=\frac{1}{2\pi i}\int_\Gamma R_t(z)udz
        $$
        for $u \in C^\infty(M,g_0)$. By \autoref{resolvent}, $F_t(\Gamma)u$ depends continuously on $t$.

        Let $\ell=\dim \mathcal F_{0}(\Gamma)$ be the multiplicity of the eigenvalue $\lambda$ and we denote the eigenfunctions of $P_{g_0}$ as $u_1(0), \dots, u_\ell(0)$ with respective eigenvalues of $\Lambda_1(0),\dots,\Lambda_\ell(0)$. Since $F_{t}(\Gamma)u$ depends continuously on $t$, for all $u\in C^\infty(M,g_0)$, we know that 
        $$
        \lim_{t \to 0}||F_{t}(u_j(0))-u_j(0)||=0
        $$
        for $j=1,\dots, \ell$. For some smaller $\delta$, $F_t(\Gamma)(u_1(0)),\dots, F_t(\Gamma)(u_\ell(0))$ are linearly independent when $|t|<\delta$, and hence 
        $$
        \dim \mathcal{F}_t(\Gamma)\geq \ell.
        $$

    Now we define $\ell'=\limsup_{t \to 0}\dim \mathcal{F}_t(\Gamma)$. Choose a sequence $\{t_i\}_i$ where $t_i \to 0$ so that $\dim F_{t_i}(\Gamma)=\ell'$, and let $m(t_i)u_{j}(t_i)$ for $j=1,\dots \ell'$ be the eigenfunctions that span $\mathcal{F}_{t_i}(\Gamma)$ with respective eigenvalues $\Lambda(t_i) \leq \cdots\leq \Lambda_{\ell'}(t_i)$. By \autoref{sobolbound},
    $$
    ||m(t_i)u_j(t_i)||_m \leq c_0 (|\Lambda_j(t_i)|+1)
    $$
    is bounded for all $j=1,\dots, \ell'$. Since the Sobolev embedding is compact for $k<m$, we may choose a subsequence so that $m(t_i)u_j(t_i)$ converges in Sobolev norm $||\cdot||_{k}$ for all $j=1,\dots, \ell'$.

    Denote $v_j$ be the limit of these eigenfunctions and observe that since $m(0)$ is the identity,
    $$
    P_{g_0}v_j=A(0)v_j =\lim_i A(t_i)m(t_i)u(t_i)=\lim_i \lambda_j(t_i)m(t_i)u(t_i).
    $$
    It follows that $v_1, \dots,v_{\ell'}$ are linearly independent eigenfunctions of $P_{g_0}$ that belong to $\mathcal{F}_0(\Gamma)$. Hence for $|t|<\delta$, the equality follows from
    $$
    \dim \mathcal{F}_0\ge \ell'
    $$
    \end{proof}
    Using the above Lemma, we prove the perturbation theory of the GJMS operators of order $N$ on the manifolds of dimension $n$.  
    \begin{proof}[Proof of \autoref{onesided}]
    We use \cite[Lemma A.4]{BU83} to verify that $A(t)$ depends $t$ analytically. Here, we want to show that for any coordinate charts covering $M$, we can write 
        $$
        A(t)=\sum_{|\alpha|\leq 2N}A_\alpha(t,x)\partial^\alpha
        $$
    where $A_\alpha(t,x)$ is analytic.

    We only need to verify that $m(t)$ is analytic since we know that $P_{g(t)}$ is analytic due to \autoref{Panalytic}. Indeed, if $|t|$ is small enough, we can write $$g^{-1}(t)=\sum_{k\ge 0}(-1)^k(g_0^{-1}G(t))^kg_0^{-1}$$ where $G(t)=\sum_{k\ge 1}g^{(k)}t^k/k!$ and $g^{(k)}=\frac{d^k}{dt^k}g(t)|_{t=0}$.

    Recall that we can write given the local chart, we can write $P_{g(t)}$ as
    $$
    P_{g(t)}=\sum_{|\alpha|\leq 2N} a_\alpha(t,x)\partial^\alpha
    $$
    where $a_\alpha$ is analytic. After arranging the terms, we have the series 
    $$
    A(t)=A^{(0)}+tA^{(1)}+t^2A^{(2)}+\dots
    $$
    that converges for $|t|<\delta$ for some $\delta$. We also know that $A_0$ is bounded below uniformly in $t$ due to the ellipticity $P_{g_0}$. 
    
    The eigenfunctions constructed by \cite[Theorem A.2]{BU83} are not necessarily real-valued. Let's denote the set of orthonormal eigenfunctions $\{\phi_1(t),\dots \phi_\ell(t)\}$ corresponding to $\{\Lambda_1(t),\dots, \Lambda_\ell(t)\}$. Suppose we have $\tilde{\ell}$ eigenvalue functions that are not repeated for $|t|<\delta$. Choose $\{\varphi_{1}(0) ,\dots, \varphi_{m_i}(0)\}$ corresponding to $\Lambda_i(0)$ with multiplicity $m_i$, which is still orthonormal set. Then we collect all the real and imaginary parts of the eigenfunctions that are linearly independent. Note that they are still eigenfunctions of $\Lambda_i(t)$ for small $\delta$.
    
    Let's use Gram-Schmidt process on this set of eigenfunctions to construct real-analytic and real-valued eigenfunctions that are orthonormal. Note that after Gram-Schmidt, the functions are still real-analytic because making $\delta$ smaller if necessary, the eigenfunctions that we chose are linearly independent. Repeat this process for any distinct eigenvalue. After this process, there are in total $\ell$ number of orthonormal eigenfunctions. Then choose $u_i=m(t)^{-1}\varphi_i$, which is the real analytic and real valued eigenfunctions.  This gives $(i)$ and $(ii)$.

    We prove $(iii)$. Note that $P_{g_0}$ has discrete eigenvalues. Hence there exists a circle around $\lambda$ that its interior and the boundary do not include any other spectrum of $P_{g_0}$ and we apply \autoref{separate} to conclude $(iii)$.
    \end{proof}

Hence, we define one-sided derivatives using the analyticity of the eigenvalues. Here, we write $E_k(g)=\ker(P_{g}-\lambda_k(g)I)$ as the eigenspace corresponding to $\lambda_k(g)$. And define a mapping $\Pi_k:L^2(M,g) \to E_k(g)$ as the orthogonal projection on $E_k(g)$. Essentially, $\Pi_k\circ P_{g}$ is a self map from $E_k(g)$ to itself. The next theorem states that under analytical perturbation of metrics, we may choose a path so that we have left and right derivatives.

\begin{proposition} \label{in}
    Suppose we have same assumptions as in \autoref{onesided}. Let $g(t)$ be an analytic perturbation of Riemannian metrics of $g_0$ where $t \in [-\delta, \delta]$. Suppose that there exists $\epsilon$ such that $\lambda_k$ is the only spectrum in $(\lambda_k-\epsilon,\lambda_k+\epsilon)$. Then $t \mapsto \lambda_k(P_{g(t)})$ admits left and right derivatives at $t=0$ and we may write 
    $$
    L(\lambda_k')=\frac{d}{dt}\lambda_k(P_{g(t)})\bigg|_{t=0^-} \text{ and } R(\lambda_k')=\frac{d}{dt}\lambda_k(P_{g(t)})\bigg|_{t=0^+}.
    $$
    Moreover, suppose either of the cases $(i)$ or $(ii)$ holds:
    \begin{enumerate}[label=(\roman*)]
     \item If $\lambda_{k}(g_0)>\lambda_{k-1}(g_0)$, then $ L(\lambda_k')$ and $R(\lambda_k')$ are the greatest and smallest element in $\Pi_k \circ P_{g_0}'$ on $E_{k}(g_0)$, respectively.
    \item If $\lambda_{k}(g_0)<\lambda_{k+1}(g_0)$, then $ L(\lambda_k')$ and $R(\lambda_k')$ are the smallest and greatest element in $\Pi_k \circ P_{g_0}'$ on $E_{k}(g_0)$, respectively.   
  \end{enumerate}
Then $L(\lambda_k')$ and $R(\lambda_k')$ are eigenvalues of the linearized operator $\Pi_k \circ P_{g_0}'$.
\end{proposition}	
\begin{proof}[Proof of \autoref{in}]
    By \autoref{onesided}, under the analytic perturbation of metrics denoted as $g(t)$ and $g(0)=g_0$, we have analytic perturbation of $P_{g(t)}$ and its corresponding eigenfunctions and eigenvalues. 
    
    We first establish the existence of the one-sided derivatives of $\lambda_k(P_{g(t)})$ at $t=0$. Consider an analytic perturbation of the background metric $g_0$ denoted by $g(t)$ for $|t|<\delta$. Let's choose $\delta>0$ so that the only  spectrum of $P_{g_0}$ lies in the the interval $(\lambda_k(P_{g_0})-\epsilon, \lambda_k(P_{g_0})+\epsilon)$ is $\lambda_k(P_{g_0})$.

Due to the fact that, the functional $t \mapsto \lambda_k(P_{g(t)})$ is continuous and $t \mapsto \Lambda_i(t)$ is analytic, we have that 
    \begin{equation} \label{partial}
        \lambda_k(P_{g(t)})=
        \begin{cases}
            \Lambda_\alpha(t) \qquad \text{if $t \in [-\delta,0)$}\\
            \Lambda_\beta(t)\qquad \text{if $t\in[0,\delta)$}
        \end{cases}
    \end{equation}
    for some $\alpha$ and $\beta$ in $\{1,\cdots,\ell\}$ and for small enough $\delta$. Here, $\alpha$ and $\beta$ are not necessarily equal since $\Lambda_i(0)=\lambda_k(P_{g_0})$ for every $i=1,\cdots, \ell$. In the case of \eqref{partial}, the left derivative at $t=0$ is $\Lambda_\alpha'(0)$ and the right derivative at $t=0$ is $\Lambda_\beta'(0)$. This establishes the existence of left and right derivatives of eigenvalue. 

    It is sufficient to prove $(i)$ since $(ii)$ follows similarly from the proof of $(i)$. 

Due to the choice of $\delta$, we have $\lambda_{k-1}(P_{g_0}) <\Lambda_i(t)$ for all $t \in (-\delta,0)$ and for all $1\leq i\leq \ell$. Moreover, $\lambda_k(P_{g(t)})=\Lambda_\alpha(t)$ for $t \in(-\delta,0]$ in \autoref{partial}, which is $\min\{\Lambda_1(t),\dots, \Lambda_\ell(t)\}$. Then the left derivative can be computed as  
$$
\frac{d}{dt}\lambda_k(P_{g(t)})\bigg|_{t=0^-}=\lim_{t \to 0^{-}}\frac{\lambda_k(P_{g(t)})-\lambda_k(P_{g_0})}{t}=\Lambda_\alpha'(0) .
$$

Similarly, for $t \in (0,\delta)$, we can assume that $\lambda_{k-1}(P_{g_0})< \Lambda_\beta (t)$. Then we have $\lambda_k(P_{g(t)})=\Lambda_\beta(t)$ for $t \in (0,\delta)$, which is $\min\{\Lambda_1(t),\dots, \Lambda_\ell(t)\}$. We calculate the right derivative as
$$
\frac{d}{dt}\lambda_k(P_{g(t)})\bigg|_{t=0^+}=\lim_{t \to 0^{+}}\frac{\lambda_k(P_{g(t)})-\lambda_k(P_{g_0})}{t}=\Lambda_\beta'(0) .
$$

Next, we show that $E_k(P_g)$ is also an eigenbasis linearized operator $P_g'$ with the eigenvalues $\{\Lambda_1'(0), \dots, \Lambda_\ell'(0)\}$. Differentiating and evaluating \eqref{varlaplace} at $t=0$, we attain
    \begin{equation}\label{pdevar}
        P_g'u_i+P_gu_i'=\Lambda_i'u_i+\lambda_k(P_g)u_i'
    \end{equation}
    where $\Lambda_i'(0)=\Lambda_i'$ and $u_i'=u_i'(0)$. Multiply $u_j$ and  subtract $u_i'P_gu_j$ to both sides, and then integrate over $M$. Then the last term in the right hand side vanishes and \autoref{pdevar} becomes
     \begin{align*}
        &\int_M \left( u_jP_g'u_i+u_jP_gu_i'-u_i'P_gu_j\right)dv_g =\Lambda_i'(0)\delta_{ij},
    \end{align*}
    where $\delta_{ij}$ is the Kronecker delta function. Here, the last two terms of the left hand side vanishes due to self-adjointness. Then we have
    \begin{equation} \label{linearized}
    \int_Mu_jP_g'u_idv_g=\Lambda_i'(0)\delta_{ij}.
    \end{equation}     
\end{proof}
If we perturb the metric within the conformal class, we can write $\Lambda'_i(0)$ explicitly. By writing $g(t)=e^{2\omega_t}g_0$, we prove later in the proof of \autoref{quadratic} (See \autoref{eigenvaluediff2}) that 
	\begin{equation} \label{eigenvaluediff1}
	\Lambda_i'(0)=-2N\lambda_k(P_g)\int_M\omega'u_i^2 dv_g,
	\end{equation}
Due to analytic perturbation, $\{\Lambda_1(t),\dots, \Lambda_\ell(t)\}$ can be well-defined and analytic for $|t|<\delta$. Note that unless $\lambda_k(P_{g(t)})$ has a critical point at $t=0$, the left and right derivatives might not be equal to each other depending on the path that we choose.

Based on our discussion in this section, we define metrics that corresponds to the critical point of the eigenvalue functional in the conformal class. We again assume that the volume is preserved under the perturbation.
\begin{defn} \label{critical}
        Suppose we have same assumptions as in \autoref{onesided}. A metric $g$ is said to be \textit{critical} for the eigenvalue functional $t \mapsto\lambda_k(P_{g(t)})$ if for any volume-preserving analytic perturbation $g(t)$ where $g_0=g(0)$, we have 
    \begin{equation} \label{eigneg}
        L(\lambda_k')\times R(\lambda_k') \leq 0.
    \end{equation}
\end{defn}
 The condition \eqref{eigneg} gives \autoref{confext}. Assume that $g_0$ is a critical metric and either of the cases holds
    \[
    L(\lambda_k')\leq 0\leq  R(\lambda_k')  \text{ or }R(\lambda_k')\leq 0\leq  L(\lambda_k')    
    \]
holds. Then the analyticity of the eigenvalue functional implies that for $|t|$ sufficiently small, we can approximate the eigenvalue respectively as follows:
	$$
    	\lambda_k(P_{g(t)}) \geq \lambda_k(P_{g_0})+o(t)\text{ or }  \lambda_k(P_{g(t)}) \leq \lambda_k(P_{g_0})+o(t).
	$$
We immediately see that if the eigenvalue is either a local minimum or local maximum at the metric $g_0$, then $g_0$ is critical.

\section{Proof of main results}\label{proofs}

	The proof of \autoref{main} is mainly from the idea by El Soufi and Ilias \cite[Theorem 2.2]{EI03}. Per\'{e}z-Ayala \cite{P22} applied this method to the Paneitz operator. The simplification in this paper relies on the fact that the bilinear operator is the remainder of the Leibniz formula of the general GJMS operators.
    
   \begin{proof}[Proof of \autoref{main}]   
    First, let's show the existence of eigenfunctions $u_1, \dots , u_\ell$ satisfying the condition $\sum_{i=1}^\ell u_i^2\equiv 1$. Recall that we define $E_k(g_e)$ as the space of eigenfunctions of $\lambda_k({g_e})$. Moreover, we assume that the eigenfunctions are normalized such that if $f \in E_k(g_e)$, then $||f||_{L^2(M,g_e)}=1$. Define a convex hull $H$ of $E_k(g_e)$ defined as
    $$
    H=\left\{\sum_{i \in \mathcal{I}}a_iu_i^2 \,\middle|\, \text{$a_i \geq 0$, $\sum_{i\in\mathcal{I}}a_i=1$, $u_i \in E_{k}(P_{g_e})$ and $|\mathcal{I}|<\infty$ }\right\},
    $$
    as a subset of $C^\infty(M)$ endowed with norm $L^2(M,g_e)$. The first step is to argue by contradiction that $\V(M,g_e)^{-1} \in H$. If $\V(M,g_e)^{-1} \in H$, then choose the eigenfunctions by scaling then we are done. Suppose not, that is, $\V(M,g_e)^{-1} \notin H$. Note that $\{\V(M,g_e)^{-1}\}$ and $H$ are disjoint convex subsets. They are also compact since convex hull of a compact set in a finite-dimensional space is compact. By the classical separation theorem based on Hahn-Banach theorem, there exists a function $\omega$ in $L^2(M,g_e)$ such that for all $f \in H$,
    \begin{equation}
     \int_M f \omega\,dv_{g_e} \leq 0<\V(M,g_e)^{-1}\int_M\omega \,dv_{g_e}, \label{sepa}
	\end{equation}
    which means that $H$ and $\{\V(M,g_e)^{-1}\}$ are separated. Without loss of generality, by the smooth approximation, we may assume that $\omega$ is smooth.

   Define an auxiliary function $\tilde{\omega}$ as
        \begin{equation}
	\tilde{\omega}=\omega-\V(M,g_e)^{-1}\int_{M}\omega dv_{g_e} \label{omegatil}
	\end{equation} and let $g(t)$ be some analytic conformal perturbation of $g_e$ given by
	\[
        g(t)=\left(\frac{\V(M,g_e)}{\V(M,e^{2t\tilde\omega}g_e)}\right)^{1/N}e^{2t\tilde{\omega}}g_e,
	\]
    where $g(t)$ preserves volume. Indeed, we have 
    $$\frac{d}{dt}\V(M,e^{2t\tilde\omega}g_e)\bigg|_{t=0}=2N\int_M\tilde \omega dv_{g_e}=0.$$ To see the sign of the derivative of the eigenvalue \eqref{eigenvaluediff1} at the critical metric, it is sufficient to check the sign of 
\begin{equation} \label{eigder}
\int_M \tilde{\omega}u^2dv_{g_e}=\int_M\omega u^2dv_{g_e}-\V(M,g_e)^{-1} \left(\int_M\omega dv_{g_e}\right) \left(\int_Mu^2dv_{g_e}\right)
\end{equation}
     Due to \eqref{sepa}, the sign of the right hand side of \eqref{eigder} is negative. This means that $ L(\lambda_k')\times R(\lambda_k') >0$, which contradicts our condition \eqref{eigneg}.
    \end{proof}
   \begin{proof}[Proof of \autoref{cormain}]
	Suppose we have a conformally extremal metric $g_e$ inside the conformal class. We prove by contradiction that the eigenspace corresponding to $g_e$ has at least dimension 2. Suppose that the eigenspace has only dimension 1. Then by \autoref{main}, we know that the only eigenfunction inside the eigenspace is constant function, which is a contradiction unless $\lambda_k(P_{g_e})=0$.
   
	Next, let's obtain the formula in \autoref{eigformula}. Suppose that 
        \begin{align*}
        &0=P_{g_e}(1)=P_{g_e}\left( \sum_{i=1}^\ell u^2_i\right)=\sum_{i=1}^\ell P_{g_e}\left(u^2_i\right) 
        \\&=\sum_{i=1}^\ell\left(2u_i^2 \lambda_k(P_{g_e})-2B_g(u_i, u_i) \right),
        \end{align*}
    where the last line is due to the equation \autoref{PandB}. Since $\sum^\ell_{i=1}u_i^2=1$, we have 
    $$
    \lambda_k(P_{g_e})=\sum_{i=1}^\ell B_g(u_i, u_i).
    $$    
    \end{proof}
        \begin{proof}[Proof of \autoref{partialsuf}]
     By normalizing the eigenfunctions if necessary, we assume that $\sum_{i=1}^\ell u_i^2=1$. Here, multiplying both sides by an arbitrary test function $f$ and integrating on $M$ with respect to $g_e$, we have 
    \begin{equation} \label{norm1}
        \sum_{i=1}^\ell\int_Mf u_i^2dv_{g_e}=\int_Mfdv_{g_e}.
    \end{equation}
    Define an analytic, volume-preserving perturbation of $g_e$ in $[g]$ as $e^{2\omega_t}g_e$. By writing $\omega'=\frac{d}{dt}\omega_t|_{t=0}$, the volume preserving property gives  
    \begin{equation} \label{0condition}
        0=\int_M\left(2N\frac{d}{dt}\omega_t\bigg|_{t=0}\right)dv_{g_e}=\int_M2N\omega'dv_{g_e}.
    \end{equation}
    Therefore, combining both \eqref{norm1} and \eqref{0condition}, we have 
    \[
    \sum_i^\ell\int_M2N\omega'u_i^2dv_{g_e}=\int_M 2N\omega'dv_{g_e}=0
    \]
    Since we assume that $\lambda_k(P_{g_e}) \neq0$, the derivative of the eigenvalue defined says that we have a conformally extremal metric as in \autoref{critical}.
    \end{proof}

    \begin{proof}[Proof of \autoref{partialsufcor}]
    Recall that $U=(u_1,\cdots,u_\ell):(M,g) \to (\mathbb{S}^{\ell-1},g_r)$ solves the equation $$P_g(U)=\mathcal E_g(U)U.$$ Assume $\mathcal E_g(U)$ is always positive on $M$. Write $g_U$ as the metric $\mathcal E_g(U)^{\frac{1}{N}}g$. Then 
    \[
        P_{g_U}u_i=u_i
    \]
    where we can see that $u_i$ of $U$ are eigenfunctions of $P_{g_U}$ with corresponding eigenvalue as 1. If 1 is the $k$-th eigenvalue, written as $\lambda_k(P_{g_U})$, we have chosen this $k$ such that $\lambda_k(P_{g_U})>\lambda_{k-1}(P_{g_U})$. Note that due to our assumption that $\sum_{i=1}^\ell u_i^2=1$, we have $P_{g_U}\left(\sum_{i=1}^\ell u_i^2\right)=0$. (Here, we used the same technique in the proof of \autoref{main}.) So we have that $\mathcal E_{g_U}(U)=\lambda_k(P_{g_U})$ is constant. And we have the conclusion from \autoref{partialsuf}.

    Now assume $\mathcal E_g(U)$ is always negative on $M$. By taking the metric $g_U=(-\mathcal E_g(U))^{\frac{1}{N}}g$, the coordinate functions of $U$ are exactly eigenfunctions of $P_{g_U}$ with corresponding eigenvalue $-1$. The remaining part follows similarly from the previous case.  
    \end{proof}
    
    Next, we investigate the necessary and sufficient condition of critical metrics when the metric varies analytically. We define the quadratic form corresponding to the linearization of $P_g$.
    \begin{lemma} \label{quadratic}
        Let $(M,g)$ be a closed $2N$-dimensional manifold endowed with a Riemannian metric $g$ and $u\in  C^\infty(M)$. Consider a family of analytic perturbation of $g$ denoted as $\{g(t):|t|<\delta\}$ with $h=\frac{d}{dt}g(t)|_{t=0}$. Then there exists an operator $P_h'=\Pi_k\circ P_g':E_k(g) \to E_k(g)$ where $ P_h'$ is symmetric in $L^2$-norm with respect to $g$ and the quadratic form $Q_h$ corresponding to the operator 
        \begin{equation} \label{quadfor}
             Q_h(u)=\int_MuP_{h}'udv_g=-
        \int_M\langle h, \tau_g(u,u) +\frac{1}{2} \lambda_ku^2 g\rangle dv_g,
        \end{equation}
        where $\tau_g$ is a symmetric traceless 2-tensor. 
    \end{lemma}

    \begin{proof}
     Recall that the linearization of $dv_g$ is written as 
    \begin{equation} \label{linmeasure}
    \frac{d}{dt} dv_{g(t)}\bigg|_{t=0}=\frac{1}{2}\langle  g, h\rangle dv_g
\end{equation} Let $u$,$w$ be arbitrary functions in $ E_k(P_g)$. We use the product rule to compute the linearization of $P_g$ using \eqref{linmeasure},
    \begin{align}
     &\int_M uP_h'w\,dv_g=\frac{d}{dt}\int_M uP_{g(t)} w dv_{g(t)}\bigg|_{t=0}-\int_M uP_{g(t)}w\frac{d}{dt}dv_{g(t)}\bigg|_{t=0} \label{linpart}
     \end{align}
    
    To compute the left hand side in \autoref{linpart}, recall that one can $P_g$ by the perturbation theory stated in \autoref{onesided} under analytic perturbation. For $u,w \in C^\infty(M,g)$, $uP_{g(t)}w$ can be universally written as
            $$
            uP_{g(t)}w=\sum_{|\alpha|\le 2N}ua_\alpha(t,x)\partial^\alpha w,
            $$
    where $a_\alpha(t,x)$ are real-analytic by \autoref{Panalytic}. So we write the linearization of $a$ as $\dot a$. Moreover, we later show in \autoref{computation} that $a$ can be written as a polynomial $F$ of $g$, $g^{-1}$ and partial derivatives of $g$. Then there exists another polynomial $\tilde F$ such that
    $$
    \dot a_\alpha(x)=\frac{d}{dt}F(g(t),g(t)^{-1},\partial^\beta g(t))\bigg|_{t=0}
    =\tilde F(h,\partial^{\beta} h ),
    $$
    where the coefficients of $\tilde F$. For fixed $\alpha$, the polynomial $\tilde F$ consists of the terms in the form of $u\tilde a(x)^{ij}\partial^\alpha v\partial^\beta h_{ij}$. Using the integration by parts, we can move $\partial^\beta$ to the front, we can change this term into $(-1)^\beta \partial^\beta(u\tilde{a}(x)^{ij}\partial^\alpha v) h_{ij}$. And the term $(-1)^\beta\partial^\beta(u\tilde{a}(x)^{ij}\partial^\alpha v)$ is a 2-tensor. Note that $\tilde a$ depends on the initial metric $g$.

    The uniqueness of this tensor follows from applying the previous argument for each of the charts and recovering the formulas on all of $M$. Since on each of the chart, the functions vanish on the boundary on some neighborhood. And the partition of unity argument establishes the uniqueness. And the symmetry follows from the self-adjointness of $P_{g(t)}$.
 
    Add $g\lambda_k(P_g)\frac{1}{2}uw$ to this tensor and then define $\tau_g(u,w)$ as  
    \begin{equation*}
        \frac{d}{dt}\int_M uP_{g(t)} w dv_{g(t)}\bigg|_{t=0}=-\int_M \langle h, \tau_g(u,w)\rangle\,dv_g, \label{linall}
    \end{equation*}
    where $\tau_g$ is a 2-tensor. Also it is traceless if we show that $\langle \tau_g(u,w),g\rangle=0$. Suppose we take the perturbation $g(t)=e^{2\omega t}g$, $h=2\omega g$. Since the first integral in the right hand side of \autoref{linpart} does not depend on the $\tau_g$ due to the conformal invariance, we have 
    $$
    0=\int_M \omega \langle \tau_g(u,w),g \rangle dv_g
    $$
    for any smooth functions $\omega$, and hence we conclude that the tensor $\tau_g$ is traceless.
    \end{proof}

As a remark, if we limit our perturbation in the conformal class, we can write the perturbation $g(t)=e^{2\omega_t}g$ for some smooth function $\omega_t$. Then we can use the bilinear operators defined in \autoref{wodz} to obtain an expression for $P_g'=\frac{d}{dt}P_{g(t)}|_{t=0}$. This allows us to use the conformal invariance property of the bilinear operators and write 
    $$
     B'_g(u,w)=-\frac{1}{2}P'_{g}(uw)+\frac{1}{2}uP'_gw+\frac{1}{2}wP'_gu.
    $$
Then we can write the first term on the right hand side of \autoref{linpart} as 
    \begin{align*}
        \frac{d}{dt}\int_M uP_{g(t)} w dv_{g(t)}\bigg|_{t=0}&=\frac{d}{dt}\int_MB_{g(t)}(u,w)dv_{g(t)}\bigg|_{t=0}=0
    \end{align*}
    where the property $B_{e^{2\omega_t}g}(u,w)=e^{-2N\omega_t}B_g(u,w)$ was used. We write \autoref{quadfor} as
	\begin{equation} \label{lindiff}
          \int_M uP_h'u\,dv_g=-2N\lambda_k(P_g)\int_M\omega'u^2 dv_g .
	\end{equation}
By combining this with \eqref{linearized}, we have 
	\begin{equation} \label{eigenvaluediff2}
	\Lambda_i'(0)=-2N\lambda_k(P_g)\int_M\omega'u_i^2 dv_g 
	\end{equation}
for each $i=1,\dots, \ell$, which was \autoref{eigenvaluediff1}.

To give a better understanding about \eqref{quadfor}, let's see what happens depending on the sign of $Q_h$. Let $g(t)$ be a conformal deformation of a Riemannian metric of $g$ on $M$. If $Q_h$ is positive on $E_k(g)$, by \eqref{eigenvaluediff1}, there exists $\delta>0$ such that $\lambda_k(g(t))$ is increasing for all $0<t<\delta$. Analogous result holds when $Q_h$ is negative on $E_k(g)$. 

As a partial converse, assume that $\lambda_k(P_{g}) \neq 0$. Suppose that $\lambda_k(P_g)> \lambda_{k-1}(P_g)$ and there exists $u \in E_k(g)$ such that $Q_h(u)<0$, then for some $\delta>0$, $\lambda_k(P_{g(t)})<\lambda_k(P_g)$ for all $0<t<\delta$. When $\lambda_k(P_g)<\lambda_{k+1}(P_g)$ and there exists $u \in E_k(g)$ such that $Q_h(u)>0$, then $\lambda_k(P_{g(t)})>\lambda_k(P_g)$ for all $0<t<\delta$. 
	\begin{proposition}\label{indefinite}
	If $g$ is a critical metric for the function $\lambda_k(P_g)$ on $M$, then for all analytic deformation denoted as $g(t)$ that preserves the volume, the quadratic form $Q_h(u)$ is indefinite on $E_k(g)$.
	\end{proposition}
	\begin{proof}
	Choose $h \in S^{2}(T^*M)$ that $\int_M\langle g,h\rangle dv_g=0$. Then the analytic perturbation defined as  $$g(t)=\left(\frac{\V(M,g)}{\V(M,g+th)}\right)^{1/N}(g+th)$$
    preserves volume. 
    
     Suppose we have a critical metric $g(0)=g$ when $t=0$. Using \autoref{in} and \autoref{quadratic}, we know that the linearized operator $P_h'$ admits at least one nonnegative and one nonpositive eigenvalues on $E_k(g)$. This means that $Q_h$ is neither positive or negative on $E_k(g)$.
	\end{proof}

The proof follows similarly from the work of \cite{EI08}, \cite{P22}. In the proof, we introduce a function space that symmetric covariant 2-tensors denoted as $S^2(T^*M)$.
\begin{proof}[Proof of \autoref{necessary2}]
    First we prove $(i)$. Consider the convex hull $K$ in $L^2(S^2(T^*M))$ defined as 
	\[	
        \left\{\tau_g(u,u)+\frac{1}{2} \lambda_k(P_g)u^2g:u \in E_k(P_g) \text{ and }||u||_{L^2(M,g)}=1\right\}
	\]
    Our goal is to show that $\lambda_k(P_g)g/2\V(M,g) \in K$. We prove by contradiction so suppose $\lambda_k(P_g)g/2\V(M,g) \notin K$. Note that both sets are disjoint convex sets and compact by a similar reasoning. By Hahn-Banach theorem, there exists $h \in L^2(S^2(T^*M))$ such that
    $$
    \lambda_k(P_g)(2\V(M,g))^{-1}\int_M\langle h,g\rangle dv_g>0
    $$
    and
    $$
    \int_M\langle h,\tau_g(u,u)+\frac{1}{2}\lambda_k(P_g)u^2g\rangle dv_g<0
    $$
    for all nonzero $u \in E_k(P_g)$. By the standard approximation theorem, we may assume that $h$ is smooth. We define a perturbation $g(t)$ that preserves volume as 
    $$
    \frac{d}{dt}g(t)\big|_{t=0}=\tilde h \qquad \text{where }\tilde{h}=h-(2N\V(M,g))^{-1}\int_M\langle h,g \rangle dv_g \cdot g,
    $$
    which is similar in \eqref{omegatil}. Since $g$ is critical, $Q_{\tilde h}$ is indefinite on $E_k(P_g)$ by \autoref{indefinite}. By \eqref{quadfor} and our assumptions, we have
    \begin{align}
        Q_{\tilde{h}}(u)
        &=-\int_M\langle\tilde h, \tau_g(u,u)+\frac{1}{2}\lambda_k(P_g)u^2g\rangle dv_g \nonumber\\
        &=-\int_M\langle h, \tau_g(u,u)+\frac{1}{2}\lambda_k(P_g)u^2g\rangle dv_g+\frac{\lambda_k(P_g)}{2N\V(M,g_e)}\left(\int_M\langle h,g\rangle dv_{g}\right) \left(\int_Mu^2dv_{g}\right)\nonumber\\
        &>0\nonumber
    \end{align}
for all $u \in E_k(P_g)$. This contradicts the assumption that $Q_{\tilde{h}}$ is indefinite.

Since $\lambda_k(P_g)g/2\V(M,g) \in K$, there exists a collection $\{u_1, \cdots, u_\ell\} \subset E_k(P_g)$ such that 
$$
\sum_{i=1}^\ell \left( \tau_g(u_i,u_i)+\frac{1}{2}\lambda_k(P_g) u_i^2 g\right)=\frac{1}{2\V(M,g)}\lambda_k(P_g)g.
$$
By taking the inner product with $g$ on both sides, we have $\sum_i\tau_g(u_i,u_i)=0$ since $\tau_g$ is traceless. Then we have \autoref{summation}.

For $(ii)$, set any arbitrary analytic $g(t)$ of $g$ that preserves volume, and set $h= \frac{d}{dt}g(t)|_{t=0} \in S^2(T^*(M))$ and assume there exists a set of eigenfunctions $\{u_1,\dots, u_\ell\}$ associated with $\lambda_k$ that satisfies \autoref{summation}. Therefore, we have 
    \begin{align*}
        \sum_{i=1}^\ell Q_h(u_i)
        &=-\sum_{i=1}^\ell \int_M\langle h,  \tau_g(u_i,u_i)+ \frac{1}{2}\lambda_k(P_g) u_i^2 g\rangle dv_g\\
        &=-\frac{1}{2}\lambda_k(P_g)\sum_{i=1}^\ell u_i^2 \int_M\langle g,h \rangle dv_g=0
    \end{align*}
    where the last line holds since the perturbation keeps the volume fixed. Hence $Q_h$ is indefinite so $g$ is a critical metric by using \autoref{in}.
\end{proof}
As a remark when $N=1$, note that since $\lambda_k(g)>0$ for all $k>1$, we do not need $\text{sign}(\lambda_k(P_g))$ in the proof. 

    Let's prove the theorem \autoref{main2}. The proof is a generalization of the result by P\'{e}rez-Ayala \cite[Theorem 1.4]{P22}. We rely on the variational characterization for the eigenvalues of $P_{g_\omega}$ is given by 
    \begin{equation} \label{ral2}
        \lambda_k(P_{g_\omega})=\inf_{B\in S_k(C^\infty(M))}\sup_{u\in B\setminus\{0\}}\frac{\int_MuP_{g_\omega}udv_{g_\omega}}{\int_Mu^2dv_{g_\omega}},
    \end{equation}
    where $S_k(C^\infty(M))$ is the space of all $k$-dimensional subspaces of smooth functions on $M$. We rely on the result by Hassannezhad's result \cite{H11} to show this. Let's first understand the concepts introduced in the paper. Let $M$ be a manifold with metric $d$ and measure $g$ that is non-atomic. Given $\kappa>1$, $\rho>0$ and a positive integer $B$, we say that $(M,d)$ satisfies \textit{$(\kappa,B,\rho)$-covering property} if each ball of radius $r$ with $0<r\leq \rho$ can be covered by at most $B$ balls of radius $r/\kappa$.
   
    \begin{lemma}[Hassannezhad \cite{H11}]\label{Hassen}
        Let $(M,g)$ be a manifold satisfying $(2,B,1)$-covering property, with $M$ compact and $\mu_g$ a non-atomic finite measure. Then there exists a constant $k_M$, depending on the covering property, such that for all $K>k_M$, there exists a family of annuli $\mathcal A=\{A_i\}_{i=1}^K$ with the properties 
         \begin{enumerate}[label=(\roman*)] 
            \item $\mu(A_i)\geq \mu(M)/8c^2K$, where $c$ is a positive constant depending on the covering property.
            \item the family $\{2A_i\}_{i=1}^K$ is disjoint.
            \item the outer radius of each member in $\mathcal A$ is less than 1.
        \end{enumerate}    
   \end{lemma}
    \begin{proof}[Proof of \autoref{main2}]
         Let $g_\omega\in[g_0]$ be arbitrary metric in the conformal class of base metric $g_0$. Let $d_{g_0}$ be the distance function induced by metric $g_0$. The metric space $(M,d_{g_0})$ satisfies the $(2,B,1)$-covering property for some $B\in \N$ depending only on the dimension of the manifold. 

        We apply the \autoref{Hassen} to metric-measure space $(M,d_{g_0},dv_{g_\omega})$. We may choose $K>k_M$ to obtain a family $A=\{A_i\}_{i=1}^B$ of annuli satisfying the properties above. And we write $$A_i=A(x_i,r_i,R_i)=\{y\in M \,|\,r_i <d_{g_0}(x,y)\leq R_i\}.$$
        Let's construct smooth functions $u_i$ on $M$ such that 
        \begin{itemize}
            \item $u_i \equiv 1$ on $A_i$ and $0\leq u_i \leq 1$ everywhere,
            \item $\text{supp}(u_i) \subset 2A_i$,
            \item $|\nabla_{g_0}^{k}u_i|=O(1/r_i^{k})$ on $A(x_i,r_i/2,r_i)$, and
            \item $|\nabla_{g_0}^{k}u_i|=O(1/R_i^{k})$ on $A(x_i,R_i/2,R_i)$,
        \end{itemize}
    for all $k=1, \dots, 2N$. To choose such $u_i$, consider smooth functions $w_i$ and $\tilde{w}_i$, with $0\leq w_i,\tilde{w}_i\leq 1$ and 
    \[
    w_i\equiv 1 \text{ on $B_{g_0}(x_i,R_i)$, }\text{supp}(w_i) \subset B_{g_0}(x_i,2R_i), |\nabla_{g_0}^{k}w_i|\leq \frac{C}{R_i^{k}}, 
    \]
    and
    \[
    \tilde{w}_i\equiv 1 \text{ on $B_{g_0}(x_i,r_i)^c$, }\text{supp}(\tilde{w}_i) \subset B_{g_0}\left(x_i,\frac{r_i}{2}\right)^c, |\nabla_{g_0}^{k}\tilde w_i|\leq \frac{C}{r_i^{k}}, 
    \]
    where the constant $C$ is a constant that depends only on the background metric $g_0$. Then one can easily verify that $u_i=w_i\tilde{w}_i$ satisfies the four properties above. 

    For any $i$, this function $u_i$ is a valid trial function of the variational characterization. Note that under the conformal invariance, we can write \eqref{ral2} as
    \begin{align}
        \lambda_k(P_{g_\omega})
        &\leq\frac{\int_Mu_iP_{g_\omega}u_idv_{g_\omega}}{\int_Mu_i^2dv_{g_\omega}}=\frac{\int_Mu_iP_{g_0}u_idv_{g_0}}{\int_Mu_i^2dv_{g_\omega}}. \label{confinv}
    \end{align}
    Here, the denominator can be bounded below by the volume of the annulus with respect to the measure $dv_{g_\omega}$. By integrating by parts, there exists $\tilde C$ depending only on the metric such that
    $$\left|\int_Mu_iP_{g_0}u_idv_{g_0}\right|\leq \tilde C.$$ 
    Using the property $(i)$ in \autoref{Hassen}, \autoref{confinv} becomes   
    \begin{align}
    \lambda_k(P_{g_\omega})\V(M,g_\omega)&\leq 8c^2K\tilde C
    =A([g_0]).\nonumber
    \end{align}
    Since $g_\omega$ is arbitrary, the result follows.
    \end{proof}
    Given some particular background metrics, this result states that it is reasonable to find upper bounds. 

    The following results are mainly from the work by Canzani \cite{C14}, and we adapt it to our case. Here, when we say a residual set, it is defined as a complement of a meager set, which is a countable union of nowhere dense sets. 
    \begin{proposition}[Theorem 2.1 in \cite{C14}]
        For $P_{N,g}$, the set of smooth functions $\omega$ on a Riemannian manifold where all the non-zero eigenvalues of $P_{e^\omega g, N}$ are simple. Then the set of functions is a residual set in $C^\infty(M)$,
    \end{proposition} 

If $g_e \in [g]$ is conformally extremal for the $k$-th eigenvalue functional and $\lambda_k(P_{g_e})$is nonzero, then by \autoref{main}, the degeneracy of the eigenvalue gives that either $\lambda_k(P_{g_e})=\lambda_{k-1}(P_{g_e})$ or $\lambda_k(P_{g_e})=\lambda_{k+1}(P_{g_e})$. The following proposition states that we can say which case holds.

    \begin{proposition}\label{nec2}
    Let $(M,g)$ be a closed $2N$-dimensional Riemannian manifold.
        \begin{enumerate}[label=(\roman*)] 
            \item If $g_e$ is a local minimizer for \eqref{evfunctionalconf} and $\lambda_k(P_{g_e}) \neq0$, then $\lambda_k(P_{g_e})=\lambda_{k-1}(P_{g_e})$
            \item If $g_e$ is a local maximizer for \eqref{evfunctionalconf} and $\lambda_k(P_{g_e}) \neq0$, then $\lambda_k(P_{g_e})=\lambda_{k+1}(P_{g_e})$
        \end{enumerate}
    \end{proposition}

    \begin{proof}
	 We prove by contradiction. Assume that $g_e$ is a local minimizer for \eqref{evfunctionalconf}, and $\lambda_k(P_{g_e})\neq\lambda_{k-1}(P_{g_e})$. Choose an analytic perturbation $g(t)=e^{2\omega_t}g_e \in [g]$ satisfying 
	   \[
           \int_M\omega'dv_{g_e}=0
	   \]
        where $\omega'=\frac{d}{dt}\omega|_{t=0}$ and this makes $g(t)$ preserves the volume. If $\ell$ denotes the multiplicity of $\lambda_k(P_{g_e})$, then by applying \autoref{onesided}, there exists a set of eigenvalue functionals written as $\Lambda_1(t),\cdots, \Lambda_\ell(t)$, and their corresponding eigenfunctions written as $u_1(t), \cdots, u_\ell(t)$ that satisfy \eqref{varlaplace} and \eqref{normalize}.

        By $(i)$ in \autoref{in}, there exists a number $\delta>0$ such that the following holds:
        \begin{equation} \label{eigall}
            \frac{d}{dt}\lambda_k(P_{g})\bigg|_{t=0}=\max_{1\leq i\leq \ell}\Lambda_i'(0)=\min_{1\leq i\leq \ell}\Lambda_i'(0)=0
        \end{equation} 
	where $|t|<\delta$. This implies that $\Lambda'_i(0)=0$ for $i=1,\cdots, \ell$. Recall that by writing $u_i=u_i(0)$ for $i=1,\cdots,\ell$, the set of eigenfunctions $\{u_1,\cdots,u_\ell\}$ forms an orthonormal basis for the $k$-th eigenspace $E_k(P_{g_e})$. Combining \autoref{eigenvaluediff1} and \eqref{eigall}, we have that 
            \[
            \int_M\omega'u^2_idv_{g_e}=0.
            \]
        By \eqref{linearized} on the proof of \autoref{in}, we obtain 
        \[
        \int_M\omega'u_iu_jdv_{g_e}=0
        \]
        whenever $i\neq j$. Therefore, we have that 
        \[
        \int_M\omega'u^2dv_{g_e}=0
        \]
       for all $u \in E_k(g_e)$. Since $\omega$ was chosen arbitrarily and is orthogonal to constant, this means that $u$ is constant for any $u \in E_k(P_g)$, which is a contradiction to $\lambda_k(P_{g_e})$ being nonzero. The proof for $(ii)$ follows similarly.
    \end{proof} 

As a remark, suppose $\text{Ker}(P_g)$ consists only of constant functions. Then we have $\lambda_2(P_g)>\lambda_1(P_g)=0$. We can verify this by contradiction. Suppose that there exists a local minimizer $g_e$ in $[g]$ for the second eigenvalue functional. Then by \autoref{nec2}, $\lambda_1(P_{g_e})=\lambda_2(P_{g_e})=0$, which is a contradiction to our assumption of $\text{Ker}(P_g)$. This is a direct adaptation of \cite[Proposition 4.2]{EI08}.

    Note that we know that $\text{Ker}(P_g)$ consists only of constant functions and the positivity of $P_g$ are conformally invariant. If the reader wants to read more work on sufficient condition of a conformal classes with these properties, see the work by Gursky \cite{G99}. See Corollary 2.1 for a similar result by El Soufi and Ilias \cite{EI03}.

	Now we prove the \autoref{first}. We skip some of the details since they were explained well in my previous paper \cite{K25}. The goal is to construct trial functions, which we denote as $u_i$ for $i=1, \dots 2N$. For a fixed $g \in [g_0]$, we use the Rayleigh quotient
	\begin{equation} \label{ral}
    	 \frac{\int_M u_i P_g u_idv_g}{\int_M u_i^2dv_{g}},
	\end{equation}
    where $u_i $ is in $H^{N}(M)$ for all $i=1,\dots, \ell$. Let $g_\omega=e^{2\omega}g_0$ be any metric in $[g_0]$ with volume equal to $\V(M,g_0)$. The proof of \autoref{first} relies on the idea by Hersch \cite{H70}, which was later adapted by Li and Yau \cite{LY82} and El Soufi and Ilias \cite{EI86}.

	\begin{proof}[Proof of \autoref{first}]
	Given a metric $g_\omega$ in the conformal class $[g_0]$ on $\Sp^{2N}$. Our goal is to find set of trial functions. Let's denote the map consisting of coordinate functions $Y=(y_1,\dots, y_{2N+1})$. Laugesen \cite{L20b} showed the existence of a center of mass of a borel measure. A point $c(g_\omega) \in \B^{2N+1}$ is the (hyperbolic) center of mass of a borel measure induced by the metric $g_\omega$ on the sphere if
	\[
	\int_{{\mathbb S}^{2N}} T_{-c(g_\omega)}(y)  \, dv_{g_\omega}=0. \label{center}
	\]

Now we substitute these functions in the variational characterization \eqref{ral} and write $c(g_\omega)=c$. Then we have
\begin{equation}
\label{Rayleigh2}
\lambda_{2}(P_{g_\omega}) \leq \frac{\int_{\Sp^{2N}}^{}(Y\circ T_{-c} )_{j}P_{\!g_\omega\,}(Y\circ T_{-c})_{j} \, dv_{g_\omega}}{\int_{\Sp^{2N}}^{}(Y\circ T_{-c})_{j}^{2} \, dv_{g_\omega}}, \qquad j=1, \dots,2N+1.
\end{equation} 

For each $j$, we multiply by the denominator on both sides of the inequality \autoref{Rayleigh2} and sum for all $j$'s. Note that $(Y\circ T_{-c} )$ is a point on the unit sphere. This implies that the sum of the denominators is equal to the volume of $M$ with the metric $g_\omega$, that is, $ \int_{\Sp^{2N}} \sum_{j} (Y\circ T_{-c})_j^{2} dv_{g_\omega} = \V (\Sp^{2N}, g_\omega)$. Recall that by our normalization, $\V (\Sp^{2N}, g_\omega)= \V (\Sp^{2N}, g_0)$. Hence, the inequality in \autoref{Rayleigh2} can be written as,
\[
\lambda_{2}(P_{g_\omega}) \V (\Sp^{2N}, g_\omega) \leq  \int_{\Sp^{2N}}^{} \sum_{j=1}^{2N+1}(Y\circ T_{-c})_{j}P_{\!g_\omega\,}(Y\circ T_{-c})_{j} \, dv_{g_\omega}.
\]
Using the conformal invariance property of $P_{g_\omega}$, the inequality becomes
\begin{align}
 \lambda_{2}(P_{g_\omega}) \, \V(\Sp^{2N},g_\omega) &\leq\sum_{j=1}^{2N+1}\int_{\Sp^{2N}}^{}(Y\circ T_{-c})_jP_{\!g_\omega\,}(Y\circ T_{-c})_j \, dv_{g_\omega} \nonumber
\\& = \sum_{j=1}^{2N+1}\int_{\Sp^{2N}}^{}(Y)_jP_{(T_{c})^{*}g_\omega\,}(Y)_j \, dv_{(T_{c})^{*}g_\omega}, \nonumber\\
&=\sum_{j=1}^{2N+1}\int_{\Sp^{2N}}^{}(Y)_jP_{g_0\,}(Y)_j \, dv_{g_0}
\end{align}
by changing $(T_{c})^{*}g_\omega$ to the round metric $g_0$ on $\Sp^{2N}$ and using conformal invariance.

Recall that by \autoref{einstein},
$$
P_{g_0}=\Pi_{1\leq j\leq N}(\Delta_{g_0}+(N+j-1)(N-j))
$$
hence 
\begin{equation}
 \lambda_{2}(\Sp^{2N},g_\omega)\,\V(\Sp^{2N},g_\omega)\le (2N)!, \nonumber
\end{equation}
which finishes the proof of \autoref{first}. In the equality case, we know that $T_{-c}$ has to be identity functions and the equality holds if the metric is isometric to the round metric.  
	\end{proof}
\section*{Appendix: Proof of \autoref{Panalytic}} \label{computation}

The following computation is to show that the coefficients of GJMS operators of even orders $N\leq n/2$ are analytic for sufficiently small $t$. This mainly comes from the fact that the coefficients of $P_g$ are universal polynomials depending only on the metric. In the work of Fefferman and Graham \cite{FG12}, \cite{FG13}, they also describe GJMS operators as ``natural" operators. Our goal in this section is to provide clear description of the dependence on the metrics. Note that we use $\Delta=d\delta$ and Fefferman and Graham uses the opposite sign. 

Given a chart, we know that the GJMS operator can be written as 
\begin{equation*}
 P_{N,g}= \sum_{|\alpha|\leq 2N} a_\alpha(x)\partial^\alpha
\end{equation*}
where $a_\alpha(x):M \to \R$ are coefficients that depend on the metric. We want to show that these coefficients are analytic over the perturbation of the metrics. 

Based on the universal polynomial defined in \cite{FG12}, we give the definition of class of assignments of the metric.
\begin{defn}
         Let $\mathcal{P}$ be the class of assignments where given a metric on a manifold $M$ and a local chart, an assignment in $\mathcal{P}$ maps the metric $g$ as
        \begin{equation*}
            g\mapsto F(g_{ij},g^{ij}, \partial^\alpha g_{ij}) 
        \end{equation*}
        where $F$ is a polynomial with rational coefficients and $F$ depends only on the dimension of the manifold. 
\end{defn}
     Here, the expression of the polynomial is consistent throughout the different metrics and coordinate charts. By definition $\mathcal{P}$ is closed under multiplication, addition and partial derivatives. For example, for Yamabe operator $P_{1}$, we write
     $$
     P_{1,g}=\Delta_g+\frac{n-2}{4(n-1)}R_g
     $$
    where $R_g$ is a scalar curvature. Locally, we can write this as
    $$
    \Delta_g=g^{ij}\partial^2_{ij}-\Gamma^k_{ij}\partial_k,
    $$
    where the Christoffel symbol can be written as 
\begin{equation}
    \Gamma_{ij}^k=\frac{1}{2}g^{kl}\left(\partial_ig_{jl}+\partial_jg_{il}-\partial_l g_{ij}\right). \label{expandchr}
\end{equation}
   See the article by Viaclovsky \cite{V07} for the formulas that we use.  and so the coefficients are in $\mathcal P$. And scalar curvature is in $\mathcal P$ similarly.
    
Using the natural tensors defined on \cite{FG12}, Juhl \cite{J13} constructed recursive formulas of the GJMS operator. Given $n$-dimensional manifold $(M,g)$ with a Riemannian metric $g$, we consider an ambient manifold $(\R \times M \times (-\epsilon,\epsilon), \tilde{g})$ with the normal form of the metric as 
$$
\tilde{g}=2\rho dt^2+2tdtd\rho+t^2g_\rho
$$
with $t \in \R,\rho \in (-\epsilon,\epsilon)$. Then we define $h_r$ as a smooth 1-parameter family of metrics on $M$ 
\begin{equation}\label{hr}
    g_\rho=h_r,\qquad h_0=g 
\end{equation}
where $\rho=-r^2/2$. Set 
\begin{equation} \label{Wnat}
    W(r)=\left(\frac{\det h_r}{\det g} \right)^{1/4}
\end{equation}
and define a 1-parameter family $\mathcal{M}(r)$ of second order differential operators on $M$ by 
\begin{equation} \label{Mr}
    \mathcal{M}(r)=\delta(h_r^{-1}d)-U(r),
\end{equation}
where 
\begin{equation} \label{Ufor}
    U(r)=\frac{[\partial_r^2-(n-1)r^{-1}\partial_r +\delta(h_r^{-1}d)]W(r)}{W(r)}
\end{equation}
acts as a zeroth order term. The generating function $\mathcal{M}(r)$ for second order differential operators $\mathcal{M}_{2N}$ on $M$ defined for $1\leq N \leq n/2$ by 
\begin{equation} \label{Mnat}
   \mathcal{M}(r)=\sum_{1\leq N\leq n/2 }\mathcal{M}_{2N}\frac{1}{(N-1)!^2} \left(\frac{r^2}{4}\right)^{N-1}.
\end{equation}
Here, $I$ is a multi-index and $n_I$ is a constant that only depends on $I$. We skip the definition since it will not be relevant to our discussion. 

Then the GJMS operator can be expressed as the linear combinations of $\mathcal{M}_{2I}$.
\begin{proposition}[{\cite[Theorem 1.1]{FG13}}]\label{Precur} For $1 \leq N\leq n/2$
\begin{equation*}
P_{N,g} = \sum_{|I|=N} n_I \mathcal{M}_{2I}.
\end{equation*}
where $I=(I_1, I_2, \dots, I_i)$ and $M_{2I}=M_{2I_1}\circ \cdots \circ M_{2I_i}$.
\end{proposition}
We use this construction to show that the coefficients of $P_{N,g}$ are natural tensors.
    \begin{lemma} \label{pconst}
        Let $g_\rho$ be a 1-parameter family of metrics defined from \autoref{hr}. Then $\partial^k_\rho g_{ij}|_{\rho=0}$ is in $\mathcal P$ for all $k$ where $1\leq k \le n/2-1$.
    \end{lemma}

\begin{proof}
    Based on \cite[Proposition 3.5]{FG12}, we show that all the local expression of the coefficients of $P_N$ is in $\mathcal P$. It suffices to show that a linear combination of contractions of Ricci curvature and covariant derivatives of Ricci curvature for the initial metric is in $\mathcal P$. For covariant derivatives of Ricci curvature, we prove by induction. For $k=0$, we have
$\Gamma^k_{ij}\in\mathcal P$ from \autoref{expandchr}, and
    \[
     R^l{}_{ijk}=\partial_i\Gamma^l_{jk}-\partial_j\Gamma^l_{ik}
  +\Gamma^l_{ia}\Gamma^a_{jk}-\Gamma^l_{ja}\Gamma^a_{ik}\in\mathcal P,
    \]
    hence the components of $\Ric$ written as $R_{jk}=R^i{}_{ijk}$ is in $\mathcal P$. For a nonnegative integer $k$, suppose $\nabla^k \Ric$ is in $\mathcal{P}$. Note that $\nabla^k \Ric$ a $(0,k+2)$-tensor and so we write the components as $(\nabla^k \Ric)_{i_1\dots i_{k+2}}$ then 
    \[
  (\nabla (\nabla^k \Ric))_{a i_1\dots i_{k+2}}=\partial_a(\nabla^k \Ric)_{i_1\dots i_{k+2}}
  -\sum_{q=1}^{k+2}\Gamma^b_{ai_q}\,
  (\nabla^k \Ric)_{i_1\dots i_{q-1}\,b\,i_{q+1}\dots i_{k+2}}
    \]
    where the formula comes from \cite[p.14]{V07}. Note that the covariant derivatives raises the second index. Recall that $\mathcal P$ is closed under addition, multiplication and partial derivatives, and so $\nabla^{k+1} \Ric$ is in $\mathcal P$.

    The terms of the contraction of Ricci curvature and the covariant derivatives of Ricci curvature are in $\mathcal P$ because the the contraction just decreases indices.
\end{proof}
    Next, we only need to solve the $h_r$ orders up to $r^{n/2-1}$ since we do not need higher orders. So it is useful to define  
    \[
    h_r=H_0^++H_1^+r^2+\dots +H_{n/2-1}^+r^{n/2-1}
    \]
    and we define $W$ and $U$ similarly. As a remark, by the work of Juhl \cite{J13}, the recursive formulas for higher order of GJMS operators can be constructed in the case of conformally flat cases.
\begin{lemma} \label{taylor2}
    Given a smooth 1-parameter family of metrics $h_r$ defined in \autoref{hr} and $W$ and $1/W$ in \autoref{Wnat}, we write
    \begin{enumerate}
        \item  $h_r^{-1}=\sum_{a=0}^{n/2-1} H_a^-r^{2a}$
        \item $W=\sum_{b=0}^{n/2-1}w_b^+r^{2b} \text{ and } 1/W=\sum_{b =0}^{n/2-1}w_b^-r^{2b}$.
    \end{enumerate}
    Then all of the coefficients are in $\mathcal{P}$ when $n$ is even.
\end{lemma}
\begin{proof}
    \textit{Step 1- coefficients of $\delta (h_r^{-1}d)$}
    
    According to our construction, $g_{\rho}=h_r$ where $\rho=-r^2/2$. Since $h_r$ has polynomial determined by the coefficients of $g_\rho$ (See Juhl \cite[p. 13]{J13}) up to the order of $r^{n-2}$, it is useful to define the expansions as 
        $$
        h_r=\sum_{0\le a \le n/2-1} H_a^+r^{2a} \qquad
        h_r^{-1}=\sum_{0 \le a \le n/2-1} H_a^-r^{2a}.
        $$
    where $H_a^+ =(-\frac{1}{2})^a\frac{1}{a!}\partial_\rho^ag(x,\rho)\big|_{\rho=0}$ and the second one is formally defined for now. By \autoref{pconst}, we have $H_a^+ \in \mathcal{P}$ and we use this to show that $H_a^- \in \mathcal P$. Note that $h_rh_r^{-1}=\text{id}$, we immediately obtain $h_r^{-1}|_{r=0}=H_0^{-1}=g^{-1}$. Similarly, define $H_a^-$ as follows:
    $$
    \sum_{0\leq k\leq a} H_k^+H_{a-k}^-=0
    $$
    for all $1\le a \le n/2-1 $. This means that we have
    \begin{equation}
        H_a^-=-g^{-1}\sum_{1\leq k \leq a}H_k^+H_{a-k}^-,
    \end{equation}
    where the uniqueness of $h_r^{-1}$ verifies the equality and hence $H_a^- \in \mathcal P$. One needs to be careful here since $H_a^{\pm}$ terms are not necessarily smooth for $n/2$-th term. This will not be a problem since for all orders $1\le N\le n/2$, we only need to use the terms $a=1,\dots ,n/2-1$. 

    Now we look at the term $\delta h_r^{-1}d$. Using our finite sum expression, we need  to apply this to $H_a^\pm$, which we write as $H_a$ to tackle both cases for now. We have 
        \begin{align}
            \delta(H_a d)&=\partial_i(H ^{ij}_a\partial_j)+\Gamma_{ik}^i(H^{kj}\partial_j) \nonumber\\
            &=H ^{ij}_a\partial_i\partial_j+(\partial_iH ^{ij}_a)\partial_j+\Gamma_{ik}^iH^{kj}\partial_j
            \label{expanddiv}
        \end{align} 
    where the coefficients lie in $\mathcal P$ since $H_a^{ij} , \Gamma_{ik}^i \in \mathcal P$ and $\mathcal P$ is closed under the multiplication, addition and partial derivatives.

    \textit{Step 2 - coefficients of $W(r)$}    
    
    Next, we show that the coefficients of $W(r)$ and $1/W(r)$ lie in $\mathcal P$. The expansion of $W$ and $1/W$ in even dimensions are defined as
    $$
    W(r)=w_0^+ +w_1^+r^2+\dots+w_{n/2-1}^+r^{n-2} \qquad 1/W(r)=w_0^- +w_1^-r^2+\dots+w_{n/2-1}
    ^-r^{n-2}.
    $$
    To compute the coefficients, we use Jacobi's formula: for a smooth family $A(s)$ of invertible matrices, we have 
    $$
    \partial_s\det A=\det A\cdot \tr(A^{-1}\partial_s A). 
    $$
    and apply this so that
    \begin{equation*}
           \partial_r W=\frac{1}{4}\tr(h^{-1}_r\partial_rh_r) W\label{Wrec} 
    \end{equation*}
    where a similar formula holds for $1/W$. From the previous step, the coefficients of $\tr(h^{-1}\partial_rh_r)$ are determined by $H_a^{\pm}$, and the coefficients of $W$ are determined recursively. From the definition $w_0^{\pm}=1$. For $1\le b\le n/2-1$, we have
    \begin{equation*}
        \begin{split}
                &w_{b}^{\pm}=\pm \frac{1}{4b}\sum_{d+c=b-1}T_dw_c^{\pm}, \\ &\text{where }T_d= \frac{\partial^{2d-1}} {\partial r^{2d-1}}\tr(h^{-1}\partial_rh_r)\bigg|_{r=0}=\sum_{a\ge 0}(d+1-a)\tr(H_a^-H_{d+1-a}^+).
        \end{split}
    \end{equation*}

\end{proof}
Note that due to the construction, the coefficients of $W$ and $1/W$ are defined by the finite terms of the coefficients of $h_r$ and $h_r^{-1}$. In particular, $w_b^{\pm}$ consists of polynomial of $H_a^+$ up to $a \leq b$. This is crucial in determining that the coefficients of $\mathcal{M}_{2N}$ are in $ \mathcal{P}$.

    \begin{proof}[Proof of \autoref{Panalytic}]
        Based on \autoref{taylor2}, we show that $U(r)$ is analytic in $r$ when we define $U(r)=\sum_{0\le e \le n/2-1} U_er^{2e}$ where $U_e$ is in $\mathcal P$ up to the order of $n/2-1$.

        To find the coefficients of $U$ defined in \autoref{Ufor}. We write 
        \begin{equation} \label{Ufirst}
             (\partial^2_r-(n-1)r^{-1}\partial_r)W =\sum_{1 \le b\le n/2-1} 2b(2b-n) w_b^+ r^{2b-2}.
        \end{equation}
        Moreover, the coefficients of $\Div_g(h_r^{-1}dW)$ are in $\mathcal P$ since the expansion is 
        \begin{equation}\label{Usecond}
            \delta (h_r^{-1}dW)=\sum_{0\le a,b\le n/2-1}\delta (H_a^-dw_b^+)r^{2(a+b)}
        \end{equation}
        and due to a similar argument in \autoref{expanddiv}, we have that the coefficients are $\mathcal P$.

        We multiply with the expansion $1/W$ to the sum of \autoref{Ufirst} and \autoref{Usecond}. Then we can write $U_{N-1}$ as follows by collecting the terms $r^{2N-2}$,
        \begin{equation*} 
            U_{N-1}=\sum_{b=1}^{N}2b(2b-n)w_b^+w_{N-b}^{-}+\sum_{a=0}^{N-2}\sum_{b=1}^{N-a-1}w_{N-1-a-b}^-\delta(H_a^-dw_b^+)
        \end{equation*}
    where the second term only consists of $H_a^-$ of indices up to $N-2$. Recall that $H_a^-$ is in $\mathcal P$ up to the order of $n/2-1$. The first term only consists of $H_a^+$ of index $n/2-1$ as well (Note that in case of the critical case, $N=n/2$ and the term involving $N$ vanishes). Hence the coefficients of $U_{N-1}$ are in $\mathcal P$. 
        
    By combining \autoref{Mnat} and \autoref{Mr}, we write the formulas for the coefficients of $M_{2N}$ 
        $$
        M_{2N}=4^{N-1}((N-1)!)^2(\delta(H_{N-1}d)-U_{N-1}).
        $$
    Given the recursive construction \autoref{Precur}, we next show that the composition of $M_{I_i}$ and $M_{I_j}$ has coefficients that are in $\mathcal P$ for distinct $I_i$ and $I_j$. Since both are linear differential operators, we can locally write
        $$
            M_{I_i}\circ M_{I_j}=\sum_{\alpha} a(x)\partial^\alpha\sum_{\beta}b(x)\partial^\beta =\sum_\alpha \sum_\beta a(x)\partial^\alpha( b(x) \partial^\beta)
        $$
    where $a,b\in \mathcal P$. Then the coefficients of $M_{I_i}\circ M_{I_j}$ are in $\mathcal P$ because $\mathcal P$ is closed under addition and multiplication. Hence the coefficients $P_{N,g(t)}$ are in $\mathcal P$.      
    \end{proof}

\section*{Acknowledgments}
I am supported by an RTG grant from the National Science Foundation ({\#}2135998). I would also like to thank professor Iosif Polterovich for suggesting this problem and his guidance. 

\section*{Acknowledgment of AI assistance}

The author explored some of the arguments interactively with Claude. The AI system assisted in developing and checking mathematical arguments and in computations. The author takes responsibility for all mathematical statements and proofs.

\bibliographystyle{plain}


\end{document}